\documentclass[12pt]{amsart}

\usepackage{amsfonts}
\usepackage{amscd}
\usepackage{amssymb}

\allowdisplaybreaks

\usepackage{enumerate,enumitem,hyperref,cleveref}

\date{\today}

\newtheorem{thm}{Theorem}[section]

\newtheorem{lem}[thm]{Lemma}

\newtheorem{prop}[thm]{Proposition}
\theoremstyle{definition}

\theoremstyle{remark}
\newtheorem{rem}[thm]{Remark}

\numberwithin{equation}{section}

\newcommand{\R}{\mathbb R}

\newcommand{\He}{\mathbb H}

\newcommand{\C}{{\mathbb C}}

\newcommand{\Fs}{\mathcal{F}^\lambda(\C^{2n})}
\newcommand{\G}{\mathcal{G}_\lambda}
\newcommand{\A}{\mathcal{A}^\lambda(\C^{2n})}

\renewcommand{\Re}{\operatorname{Re}}
\renewcommand{\Im}{\operatorname{Im}}
\newcommand{\tr}{\operatorname{tr}}

\usepackage{color}

\title[Representations of  Heisenberg groups]
{On the boundedness of certain operators of \\convolution type
on  Fock spaces}

\author[ S. Thangavelu]{ Sundaram Thangavelu}

\address[S. Thangavelu]{Department of Mathematics, Indian Institute of Science, Bangalore--560012, India.}
\email{veluma@iisc.ac.in}

\begin{document}
\maketitle


\begin{abstract} We find  necessary and  sufficient conditions on the convolution kernels $ \varphi $ so that certain operators on twisted Fock spaces $\Fs$ are bounded.
 \end{abstract}

\section{Introduction} \label{Sec-intro}
In this article we are concerned with the boundedness of certain operators of convolution type on the Fock spaces $\mathcal{F}^\lambda(\C^{2n})$ indexed by $ \lambda \in \R.$ 
These spaces are  weighted Bergman spaces consisting of entire functions $ F(\zeta)$ on $ \C^{2n}$ that are square integrable with respect to certain positive weight functions $ w_\lambda(\zeta).$ Given $F,  \varphi \in \Fs$ we define the convolution
\begin{equation}
 F\ast_\lambda \varphi(\zeta) = \int_{\C^{2n}} F(\zeta^\prime) \varphi(\zeta-\overline{\zeta^\prime})\,  \overline{K_\lambda(\zeta, \zeta^\prime)} \,w_\lambda(\zeta^\prime) d\zeta^\prime 
 \end{equation}
where $ K_\lambda(\zeta,\zeta^\prime) $ is the reproducing kernel associated to $ \Fs.$  The above integral need not converge for all  $F,  \varphi \in \Fs$ but when  $ \varphi $ is fixed, the convolution makes sense  for all $ F $ from a dense class of functions in $ \Fs.$ Thus we have a densely defined linear operator $ S_\varphi^\lambda F = F \ast_\lambda \varphi$ whose boundedness on $ \Fs$ received considerable attention recently, see the following papers \cite{CLSWY, CHLS, BD, GT}.\\

In these works  the authors have identified a subspace $\A \subset \Fs$ with the property that $ S_\varphi^\lambda,$ initially defined on a dense subspace of $ \Fs $ extends to the whole space as a bounded linear operator if and only if $ \varphi \in \A.$ These subspaces turn out to be the image of certain Banach algebras $ \bf A $ under a map $ G_\lambda.$ Thus $ S_\varphi^\lambda $ is bounded on $ \Fs$ if and only if $ \varphi = G_\lambda(\bf a) $ for some $ \bf a \in \bf A.$ It has also been shown that there is a unitary map $ \G$ setting up a one to one correspondence between  a Hilbert space $ \mathcal{H}_\lambda $ and $ \Fs.$ This map has the property that $ \G^\ast(F \ast_\lambda \varphi) = \G^\ast(F)\, \bf a $ if $ \varphi = G_\lambda(\bf a)$ which justifies the  definition of the convolution.\\

Though this characterisation is interesting in itself, there is no easy way to determine when a given $ \varphi \in \Fs $ belongs to the subspace $\A.$  It is desirable to have a necessary and sufficient condition on $ \varphi $ which are easy to verify, e.g. certain pointwise estimates/integrability conditions. In this article we address this problem and come up with a partial solution. More precisely, for each $ 0 < t  \leq 1/2 $ we introduce the subspaces 
\begin{equation}  
\mathcal{A}_{1/2}^\lambda(C^{2n}) \subset \A \subset \mathcal{A}_t^\lambda(C^{2n}) 
\end{equation} 
leading to the  necessary condition, $ \varphi \in \mathcal{A}_t^\lambda(C^{2n})$ for all $ 0< t< 1/2$  and a slightly stronger sufficient condition, $ \varphi \in \mathcal{A}_{1/2}^\lambda(\C^{2n})$ 
for the operator $ S_\varphi^\lambda $ to be bounded on $ \Fs.$ All these spaces $\mathcal{A}_t^\lambda(C^{2n})$ are weighted Bergman spaces and hence their membership  is defined in terms of integrability conditions/pointwise estimates.\\

We now give some details. Let us consider the case $ \lambda = 0 $ first which corresponds to analysis on the classical Fock space $ \mathcal{F}^0(\C^{2n}).$  These spaces are defined for any dimension and hence we consider $ \mathcal{F}^0(\C^n) = \mathcal{F}(\C^n) $ consisting of entire functions $ F $ that are square integrable with respect to  the Gaussian measure $ w_0(z) =  c_0\, e^{-\frac{1}{2}|z|^2}.$ For this case, we delete all subscripts and superscripts, e.g we simply write $ S_\varphi $ instead of $ S_\varphi^0 $ and so on. Boundedness of the operators 
\begin{equation}
 S_\varphi(z) = \int_{\C^n} F(w) \varphi(z-\bar{w})\, e^{\frac{1}{2} z \cdot \bar{w}} \, w_0(w)\, dw 
 \end{equation}
have been studied in \cite{CLSWY}. Let  $ d\gamma(\xi) = e^{-|\xi|^2} d\xi$ and consider the Gauss-Bargmann transform
\begin{equation}\label{gb-abel}
G: L^2(\R^n, d\gamma) \rightarrow \mathcal{F}(\C^n),\,\, \,  Gf(z) = e^{\frac{1}{4}z^2} \int_{\R^n} f(\xi)\, e^{i z \cdot \xi}\, d\gamma
\end{equation}
which is unitary. By taking $ {\bf A} =L^\infty(\R^n) $ and $ \mathcal{A}(\C^n) = G(\bf A)$ it has been proved that $ S_\varphi$ is bounded on $ \mathcal{F}(\C^n) $ if and only if $ \varphi= G(m)$ for some $ m \in L^\infty(\R^n).$\\

It is easy to see that the space $ \mathcal{A}(\C^n)$ is invariant under all  translations  in the imaginary direction. That is for any $ \varphi \in \mathcal{A}(\C^n) $ and $ a \in \R^n,$ the function $ \tau(ia)\varphi(z) = \varphi(z+ia) $ belongs to $ \mathcal{A}(\C^n).$ This leads a very simple necessary condition for $ S_\varphi $ to be bounded. A slight strengthening of the same condition also turns out to be sufficient.\\
\begin{thm}\label{thm-one} Suppose $ \varphi \in \mathcal{F}(\C^n) $ is such that $ S_\varphi $ is bounded. Then  $ \varphi$ satisfies
$$  (i) \,\sup_{ a \in \R^n}  \int_{\C^n} |\varphi(z+ia)|^2 e^{-\frac{1}{2}|z|^2} dz < \infty,  \,\,(ii)\,  \sup_{ y \in \R^n} |\varphi(x+iy)| \, e^{-\frac{1}{4}|x|^2} < \infty .$$
Conversely, if we assume that $ \varphi $ satisfies (i) and the slightly stronger condition
$$  \sup_{ y \in \R^n} \int_{\R^n}  |\varphi(x+iy)| e^{-\frac{1}{4}|x|^2} dx  < \infty ,$$
then $ \varphi = Gm $ for some $ m \in L^\infty(\R^n) $ and consequently $ S_\varphi $ is bounded.
\end{thm}

Another necessary condition can be obtained by considering the weighted Bergman space $ \mathcal{A}_t(\C^n),\,  t  \geq 0 $ defined by weight
$$ w_t(x,y) = e^{-\frac{1}{2}(|x|^2+|y|^2)} \, e^{\frac{2t}{(1+2t)}|y|^2}.$$
By considering a variant of the Gauss-Bargmann transform we can show that $ \mathcal{A}(\C^n) \subset \mathcal{A}_t(\C^n)$ for $ 0 \leq t <1/2$ leading to the following necessary condition.

\begin{thm}\label{thm-two} Given $ \varphi \in \mathcal{F}(\C^n),$ a necessary condition for $ S_\varphi $ to be bounded on  $ \mathcal{F}(\C^n) $ is that for $ 0 \leq t < 1/2$
$$ \int_{\R^{n}}\int_{\R^{n}} |\varphi(z)|^2w_t(x,y)\, dx\, dy  < \infty. $$
\end{thm} 

Let us consider the Hilbert space $ \mathcal{H} $ consisting of all tempered distributions $ f $ for which  $ f \ast q_{1/2} \in L^2(\R^n) $  where 
$$ q_t(x) = (4\pi t)^{-n/2} \,e^{-\frac{1}{4t}|x|^2} $$
is the heat kernel associated to the standard Laplacian $ \Delta $ on $ \R^n.$  We equip $ \mathcal{H} $ with the obvious norm : $ \|f\|  = \|f \ast q_{1/2}\|_2.$ Note that  the Fourier transform $ f \rightarrow \widehat{f} $  embeds  $ \mathcal{H} $ isometrically into a subspace of $ L^2(\R^n,\, d\gamma).$   For $ f \in \mathcal{H} $ we define $ \mathcal{G}(f) = G(\widehat{f}) $ so that the image of $ \mathcal{H} $ under $ \mathcal{G}$ is a subspace $ \mathcal{F}_0(\C^n) $ of the Fock space $ \mathcal{F}(\C^n).$  Instead of $ G $ we  work with $ \mathcal{G} $  and note that $\varphi =\mathcal{G}(f) $ defines a bounded operator $ S_\varphi $ if and only $ \widehat{f} \in L^\infty(\R^n).$ The usefulness of $ \mathcal{G} $ stems from the fact that it allows us to define the spaces $ \mathcal{A}_t(\C^n) $ mentioned above, see Section 2.\\

We now consider the case $ \lambda \neq 0.$ The spaces $ \Fs,$ known as twisted Fock spaces, are weighted Bergman spaces corresponding to the weight function
$$ w_\lambda(\zeta) = c_\lambda e^{-\frac{1}{2} \lambda (\coth \lambda)|\zeta|^2}\, e^{\lambda [ \Re \zeta,\,\Im \zeta ]} $$
where $ [ \cdot,\cdot ]$ is the symplectic form on $ \R^{2n} $ defined by $ [ (x,u), (y,v) ] = (u \cdot y-v \cdot x).$ The reproducing kernel for this space is given by
$$ K_\lambda(\zeta,\zeta^\prime)  = e^{\frac{1}{2}\lambda (\coth \lambda) \zeta^\prime \cdot \overline{\zeta}}\, e^{-\frac{i}{2}\lambda [ \zeta^\prime,  \overline{\zeta}]} $$
where $ [(z,w),(a,b)] = (w\cdot a -z \cdot b) $ is the natural extension of the symplectic form to $ \C^{2n}.$ On $ \Fs$ we consider operators of the form
\begin{align} \label{def:convolution-operator-tiwsted-fock}
S_\varphi^\lambda F(\zeta) = \int_{\C^{2n}} F(\zeta^\prime) \, \varphi(\zeta -\overline{\zeta^\prime}) \, e^{\frac{1}{2}\lambda (\coth \lambda)(\zeta \cdot \overline{\zeta^\prime})} \,\, e^{-\frac{i}{2} \lambda [\zeta, \overline{\zeta^\prime}]} \,w_\lambda(\zeta^\prime) \, d\zeta^\prime . 
\end{align}
In a recent joint work with Rahul Garg \cite{GT} we  have obtained a necessary and sufficient condition on $ \varphi $ so that $ S_\varphi^\lambda $ is bounded. The condition is in terms of a non-commutative analogue $ G_\lambda$ of the Gauss-Bargmann transform $G.$\\

In order to describe the known results we need to set up some notation. More details will be given in later sections. We  let $ {\bf A }=B(L^2(\R^n)),$ the Banach algebra of all bounded linear operators on $ L^2(\R^n).$  We let $ e^{-tH(\lambda)}, t >0$ stand for the  Hermite semigroup generated by the scaled Hermite operator $ H(\lambda) = -\Delta+\lambda^2 |x|^2 $ on $ \R^n.$ For $ (x,u) \in \R^{2n},$ let $ \pi_\lambda(x,u) $ be the Schr\"odinger representation of the Heisenberg group $\He^n$ realised on $ L^2(\R^n).$ Then it is known that $ \pi_\lambda(z,w) $ can be defined as densely defined operators for $ (z,w) \in \C^{2n}.$ Finally we let
\begin{equation} 
p_t^\lambda(x,u) = c_n \lambda^n (\sinh t\lambda)^{-n} e^{-\frac{1}{4}\lambda (\coth t\lambda)(|x|^2+|u|^2)}
 \end{equation}
which also has a holomorphic extension $ p_t^\lambda(z,w) $ to $ \C^{2n}.$ With these notations we define the non-commutative Gauss-Bargmann transform 
$ G_\lambda: {\bf A} \rightarrow  \mathcal{F}^\lambda(\C^{2n}) $ by
\begin{align} \label{twisted-GB}
G_\lambda(M)(z,w) = p_1^\lambda(z,w)^{-1}\, \tr\left( \pi_\lambda(-z,-w) e^{-\frac{1}{2} H(\lambda)} M e^{-\frac{1}{2}H(\lambda)} \right) .
\end{align}
Let $ \A $ denote the image of $ B(L^2(\R^n)) $ under the Gauss-Bargmann transform. In  \cite{GT}, we have proved that  the operator $ S_\varphi^\lambda  $ is bounded on $ \mathcal{F}^\lambda(\C^{2n}) $ if and only if $ \varphi \in  \mathcal{A}^\lambda(\C^{2n}) .$
 In other words, $ S_\varphi^\lambda $ is bounded if and only if $ \varphi = G_\lambda(M)$ for some $M \in B(L^2(\R^n)).$ 
  However, given $ \varphi \in  \mathcal{F}^\lambda(\C^{2n})$ it is not easy to determine when it belongs to $  \mathcal{A}^\lambda(\C^{2n}).$ \\
  
  For each $ 0 < t \leq 1/2,$ let the weight function $ w_t^\lambda(\xi,\eta) $ be defined on $ \R^{2n} \times \R^{2n} $ by
 \begin{equation}\label{weight}
   w_t^\lambda(\xi,\eta) = e^{-\frac{\lambda}{2}(\coth \lambda)(|\xi|^2+|\eta|^2)} e^{\lambda (\tanh 2t\lambda)|\xi|^2} \,  e^{\lambda(\coth \lambda- \coth 2t \lambda)|\eta|^2}.
   \end{equation}
   We then define $ \mathcal{A}_t^\lambda(\C^{2n}) $ as the weighted Bergman space corresponding to $ w_t^\lambda(\xi,\eta).$ In other words, $ \mathcal{A}_t^\lambda(\C^{2n}) $ is the Hilbert space of all entire functions $ F $ on $ \C^{2n} $ for which
$$ \int_{\R^{2n}}\int_{\R^{2n}}  |F(\xi+i\eta)|^2 \, w_t^\lambda(\xi,\eta) \, d\xi \, d\eta < \infty.$$
Then it can be shown that  (see Section 3.5)
$\mathcal{A}_{1/2}^\lambda(\C^{2n}) \subset \A \subset \mathcal{A}_t^\lambda(C^{2n}) $ for $ 0 < t  <1/2.$  In view of the result of \cite{GT} mentioned above this gives a necessary condition and also a sufficient condition on $ \varphi $ so that $ S_\varphi^\lambda$ is bounded.\\

\begin{thm} Given $ \varphi \in \Fs,$ a necessary condition for $ S_\varphi^\lambda $ to be bounded on $ \Fs $ is that for all $ 0 < t < 1/2$
$$ \int_{\R^{2n}}\int_{\R^{2n}}  |\varphi(\xi+i\eta)|^2 \, w_t^\lambda(\xi,\eta) \, d\xi \, d\eta < \infty .$$
Conversely, if the above integral is finite for $ t = 1/2,$ then $ S_\varphi^\lambda $ is bounded on $ \Fs .$
\end{thm}

This theorem will be proved in Section 3.5.The proof depends on a curious relation between Hermite and special Hermite semigroups and the Gutzmer's formula for the Hermite expansions. In the course of the proof we will show that $\mathcal{A}_{1/2}^\lambda(\C^{2n})$ is precisely the image of $ \mathcal{S}_2, $ the space of Hilbert-Schmidt operators on $ L^2(\R^n),$ under the Gauss-Bargmann transform $G_\lambda.$\\

 Thus if $ \varphi \in \mathcal{A}_{1/2}^\lambda(\C^{2n})$ then there exists a Hilbert-Schmidt operator $ T $ on $ L^2(\R^n)$ so that
$$\varphi(z,w) = p_1^\lambda(z,w)^{-1} \tr \left(\pi_\lambda(-z,-w)\, e^{-\frac{1}{2}H(\lambda)} T e^{-\frac{1}{2}H(\lambda)} \right).$$
As $ T = \pi_\lambda(f) $ for a unique $ f \in L^2(\R^{2n}),$ by appealing to  the inversion formula for the Weyl transform we get
$$\varphi(z,w) = p_1^\lambda(z,w)^{-1} \, p_{1/2}^\lambda \ast_\lambda f \ast_\lambda p_{1/2}^\lambda(z,w) .$$
This suggests that we consider the semigroup  $ T_t f (x,u)= p_{t}^\lambda \ast_\lambda f \ast_\lambda p_{t}^\lambda(x,u) $ which turns out to be none other than the Hermite semigroup generated by $ H(\lambda) $ on $ \R^{2n}.$ Thus the function $ p_1^\lambda\, \varphi $ belongs to the image  of $ L^2(\R^{2n}) $ under the Segal-Bargmann transform  taking $ f $ into the holomorphic function $ e^{-\frac{1}{2}H(\lambda)}f(z,w).$ The image of $ L^2(\R^{2n}) $ under the semigroup $ e^{-tH(1)} $ is known to be a weighted Bergman space corresponding to the weight function 
$$ U_{t}(z,w) = c_n (\sinh 4t)^{-2n} \, e^{- (\coth 2t) (|y|^2+|v|^2)+ (\tanh 2t)(|x|^2+|u|^2)} .$$
Therefore, it follows that $ \varphi $ is square integrable with respect to $ |p_1^\lambda(z,w) |^2\, U_{\lambda /2}(\sqrt{\lambda}(z,w)) $ 
which after simplification reduces to $ w_{1/2}^\lambda(\xi,\eta)$ proving the sufficiency part of the theorem.\\

We conclude this introduction with a brief description of the plan of this paper. In the next section we deal with operators $ S_\varphi$ on the classical Fock space $ \mathcal{F}(\C^n) $ which corresponds to the case $ \lambda =0.$ As one can see the proofs are quite easy in this case. In Section 3 we consider the operators $ S_\varphi^\lambda $ on $ \Fs$ in which case the proof requires quite a bit of preparation. After recalling the required preliminary material the main result is proved in Section 3.5.

\section{The operators $ S_\varphi$ on the Fock spaces $ \mathcal{F}(\C^n)$}

In this section we prove the results  on  operators of the form $ S_\varphi$ stated in the introduction.  We recall that these operators studied by several authors \cite{CLSWY,CHLS,BD}  are defined by
$$  S_\varphi F(z) = \int_{\C^n} F(w)\, \varphi(z-\bar{w})\, e^{\frac{1}{2} z\cdot \bar{w}}\, e^{-\frac{1}{2}|w|^2}\, dw$$
where $ \varphi \in \mathcal{F}(\C^n).$ As stated in the introduction,  $ S_\varphi $ is bounded on the Fock space if and only if $ \varphi = Gm$ for some $ m \in L^\infty(\R^n).$ Here $ G $ is the Gauss-Bargmann transform defined by the equation  \eqref{gb-abel}. The proofs of Theorems \ref{thm-one} and \ref{thm-two} are rather easy which we present now.\\

{\bf Proof of Theorem \ref{thm-one}} As mentioned above,  the operator  $ S_\varphi $ is bounded on $ \mathcal{F}(\C^n) $ if and only if
$$ \varphi(z) = e^{\frac{1}{4}z^2} \int_{\R^n}  m(\xi) \, e^{-|\xi|^2}\, e^{iz \cdot \xi}\, d\xi$$
for some $ m \in L^\infty(\R^n).$
A simple calculation shows that 
$$ \varphi(z+ia) = e^{\frac{1}{4}z^2} \int_{\R^n}  m(\xi- \frac{1}{2}a) \, e^{-|\xi|^2}\, e^{iz \cdot \xi}\, d\xi $$
for any $ a \in \R^n.$ In particular, 
$$ \varphi(x+ia) = e^{\frac{1}{4}|x|^2} \int_{\R^n}  m(\xi- \frac{1}{2}a) \, e^{-|\xi|^2}\, e^{ix \cdot \xi}\, d\xi.$$
As $ m $ is bounded we get  two necessary conditions on $ \varphi :$  for  any $ a \in \R^n $
$$  \int_{\C^n} |\varphi(z+ia)|^2 e^{-\frac{1}{2}|z|^2} dz \leq C,  \,\, |\varphi(x+ia)|\, e^{-\frac{1}{4}|x|^2} \leq C $$
where $ C $ is independent  of $ a.$
Conversely, if we assume that $ \varphi $ satisfies the first of the above conditions, then for some $ m_a \in L^2(\R^n, d\gamma) $ we have
$$ \varphi(z+ia) =  e^{\frac{1}{4}z^2} \int_{\R^n}  m_a(\xi) \, e^{-|\xi|^2}\, e^{iz \cdot \xi}\, d\xi.$$
Since  $ \varphi(z) = Gm(z),$ we have the equality
$$  e^{\frac{1}{4}x^2} \int_{\R^n}  m_a(\xi) \, e^{-|\xi|^2}\, e^{ix \cdot \xi}\, d\xi = e^{\frac{1}{4}(x+ia)^2} \int_{\R^n}  m(\xi) \, e^{-|\xi|^2}\, e^{i(x+ia) \cdot \xi}\, d\xi $$
which after  a change of variables in the right hand side integral  yields
$$   \int_{\R^n}  m_a(\xi) \, e^{-|\xi|^2}\, e^{ix \cdot \xi}\, d\xi =  \int_{\R^n}  m(\xi - \frac{1}{2}a) \, e^{-|\xi|^2}\, e^{i x \cdot \xi}\, d\xi .$$
Hence it follows that $ m_a(\xi) = m(\xi-\frac{1}{2}a) $ and  consequently $ \varphi $ satisfies the equation
$$ \varphi(x+ia) = e^{\frac{1}{4}|x|^2} \int_{\R^n} m(\xi-\frac{1}{2}a) e^{-|\xi|^2}\, e^{i x \cdot \xi}\, d\xi.$$
If we impose the stronger condition on $ \varphi ,$ namely
$$  \int_{\R^n} |\varphi(x+ia)| e^{-\frac{1}{2}|x|^2} dx \leq C $$
for all $ a \in \R^n,$ then by Fourier inversion formula
$$  m(\xi- \frac{1}{2}a) \,e^{-|\xi|^2}=  \int_{\R^n} \varphi(x+ia) e^{-\frac{1}{4}|x|^2} \, \, e^{-ix \cdot \xi}\, dx .$$
Therefore, by choosing $ \xi =0 $ we can conclude that $ m $ is bounded proving the theorem.\\

\begin{rem}  Note that $ L^\infty(\R^n) $ is a translation invariant subspace of  $ L^2(\R^n, d\gamma) $ whose image under $ G $ is $ \mathcal{A}(\C^n).$ Let $ L_0^2(\R^n, d\gamma) $ be any  translation invariant subspace which contains $ L^\infty(\R^n).$ As shown above, if $ \varphi = Gm$ with $ m \in L_0^2(\R^n, d\gamma),$ then $ \varphi_a(z) = \varphi(z+ia) $ corresponds to $ m_a(\xi) = m(\xi-a/2).$
\end{rem}
\vskip0.05in

\begin{rem} Since $ G: L^2(\R^n, d\gamma) \rightarrow \mathcal{F}(\C^n) $ is unitary we can define an action of $ \C^n $ on $ L^2(\R^n, d\gamma)$ by transferring the unitary operator $ \rho(w) $ which acts on the Fock space as
$$ \rho(w)\varphi(z) = e^{-\frac{1}{4}|w|^2}\, \varphi(z+w) \, e^{-\frac{1}{2} z\cdot \bar{w}}.$$
A simple calculation shows that
$$ \rho(u+iv)\varphi(z) =  e^{-\frac{1}{2}(|v|^2-i u \cdot v)} \, e^{iz \cdot v}e^{\frac{1}{4}z^2} \int_{\R^n}  m(\xi)\, e^{i w\cdot \xi} \, e^{-|\xi|^2}\, e^{iz \cdot \xi}\, d\xi.$$
Making the change of variables $ \xi \rightarrow \xi-v $ in the above integral and simplifying
$$ \rho(u+iv)\varphi(z) =  e^{-\frac{1}{2} |v|^2} \, e^{\frac{1}{4}z^2} \int_{\R^n}  m(\xi-v)\, e^{i \bar{w}\cdot \xi} \, e^{-|\xi|^2}\, e^{iz \cdot \xi}\, d\xi.$$
Thus  $ \rho(w) G(m) = G(\pi(w)m)$ where
$ \pi(w)m(\xi) =  e^{-\frac{1}{2}|v|^2}\, m(\xi-v)\, e^{i \xi \cdot \bar{w}}.$
\end{rem}
\vskip0.05in

We can get a different  necessary condition for $ S_\varphi $ to be bounded by looking at  another version of  the Gauss-Bargmann transform. We let $ \mathcal{H} $ stand for the space of all tempered distributions $ f $ such that $ f \ast q_{1/2} $ is in $ L^2(\R^n)$ equipped with the norm $ \| f\| = \| f \ast q_{1/2}\|_2.$ Then it is clear that $ f \in \mathcal{H} $ if and only if $\hat{f} \in L^2(\R^n, d\gamma).$ Hence $ f \rightarrow G(\hat{f}) $ is a unitary map   between $ \mathcal{H} $ and $ \mathcal{F}(\C^n).$ If $ \varphi = G(\hat{f}) $ and $S_\varphi $  is bounded, then $ \hat{f} \in L^\infty(\R^n).$  We can now prove Theorem \ref{thm-two}.\\

{\bf Proof of Theorem \ref{thm-two}} If $ \varphi $ is as above, then from the definition of $ G $ we see that
$$ \varphi(z) e^{-\frac{1}{4}z^2} = \int_{\R^n} \hat{f}(\xi) e^{-|\xi|^2} e^{i z \cdot \xi} \, d\xi.$$
By Plancherel, this immediately gives us the identity 
$$  \int_{\R^n} |\varphi(x+iy)|^2 e^{-\frac{1}{2}(|x|^2- |y|^2)} dx = c_n\, \int_{\R^n} |\hat{f}(\xi)|^2 e^{-2|\xi|^2} e^{-2 y \cdot \xi} d\xi .$$
Integrating the above with respect to $ q_{t/2}(y) dy $  we obtain
$$ (2\pi t)^{-n/2}\,\int_{\C^{n}} |\varphi(z)|^2 e^{-\frac{1}{2}|z|^2}  e^{-\frac{(1-2t)}{2t}|y|^2} dx dy = c_n\,  \int_{\R^{2n}} |\hat{f}(\xi)|^2 e^{-2(1-t)|\xi|^2} \,d\xi .$$
Under the assumption that $ \varphi $ is bounded, we have $ \hat{f} \in L^\infty(\R^n) $ and hence from the above  we get the necessary condition
$$ \int_{\C^{n}} |\varphi(z)|^2 e^{-\frac{1}{2}|z|^2}  e^{-\frac{(1-2t)}{2t}|y|^2} dx dy  \leq C_t $$
for all $ 0 < t <1.$ Replacing $ 2t $ by $ 1+2t$ we can write the above as
$$ \int_{\C^{n}} |\varphi(z)|^2 e^{-\frac{1}{2}|z|^2}  e^{\frac{2t}{1+2t} |y|^2} dx dy  \leq C_t $$
for $ 0 \leq  t < 1/2.$ This completes the proof of the theorem.

 \section{Operators $ S_\varphi^\lambda $ on the twisted Fock spaces $ \Fs.$}
  In this section we consider the operators $ S_\varphi^\lambda $ on $\Fs$ defined by the convolution
 \begin{align} \label{op-twist}
S_\varphi^\lambda F(\zeta) = \int_{\C^{2n}} F(\zeta^\prime) \, \varphi(\zeta -\overline{\zeta^\prime}) \, e^{\frac{1}{2}\lambda (\coth \lambda)(\zeta \cdot \overline{\zeta^\prime})} \,\, e^{-\frac{i}{2} \lambda [\zeta, \overline{\zeta^\prime}]} \,w_\lambda(\zeta^\prime) \, d\zeta^\prime . 
\end{align}
With $ G_\lambda $ defined in \eqref{twisted-GB} a necessary and sufficient condition on $ \varphi $ so that $ S_\varphi^\lambda $ is bounded on $ \Fs $ is that $ \varphi = G_\lambda(M) $ for some bounded linear operator $ M $ on $ L^2(\R^n).$ We are interested in finding certain integrability conditions  on $ \varphi $ which does not involve the operator $ M.$ In order to do so we make use of the connection between Hermite and special Hermite operators. Recall that  the definition of $ G_\lambda $ involves the Hermite as well as the special Hermite semigroup. We will be very sketchy in recalling the basic facts about these semigroups referring to the standard literature for more details, see e.g. \cite{T1,T2}. In what follows we assume $ \lambda \neq 0.$

\subsection{Heisenberg group and the Schr\"odinger representations} 
We begin by recalling some well known facts about the Heisenberg groups $ \He^n $ and their representations. The standard references for this section are \cite{GBF,T2}. The group $  \He^n$ is  just $ \R^n \times \R^n \times \R $ (or equivalently $ \C^n \times \R $) equipped with the group law
$$ (x,y,t)(u,v,s) = (x+u, y+v, t+s+\frac{1}{2}(u \cdot y-x \cdot v)).$$
For each $ \lambda \in \R, \lambda \neq 0,$ there is an irreducible unitary representation $ \pi_\lambda $ of $ \He^n$ on $ L^2(\R^n) $ explicitly given by 
$$ \pi_\lambda(x,y,t)f(\xi) = e^{i\lambda t} e^{i\lambda(x \cdot \xi +\frac{1}{2} x\cdot y)} f(\xi+y),\,\,\, f \in L^2(\R^n).$$
By setting $ \pi_\lambda(x,y) = \pi_\lambda(x,y,0) $ we observe that the irreducibility of $ \pi_\lambda $ has the following consequence: if a bounded linear operator $ T $ on $ L^2(\R^n) $ commutes with $ \pi_\lambda(x,y)$ for all $ x,y \in \R^n$ then $ T = cI $ for a constant $ c.$ However, there are non-trivial operators that commute either with  the family  $ \pi_\lambda(x,0) $ or with $ \pi_\lambda(0,y).$\\

The Fourier transform of a function $ f \in L^1 \cap L^2(\He^n) $ is defined to be the operator valued function
$$ \widehat{f}(\lambda) = \int_{\He^n} f(x,y,t)\, \pi_\lambda(x,y,t)\, dx\, dy\, dt = \int_{\R^{2n}} f^\lambda(x,y)\, \pi_\lambda(x,y)\, dx\, dy$$
where $ f^\lambda(x,y) $ is the inverse Fourier transform of $ f(x,y,t) $ in the central variable:
$$ f^\lambda(x,y) = \int_{-\infty}^\infty f(x,y,t)\, e^{i\lambda t}\, dt.$$
Thus we see that $ \widehat{f}(\lambda) = \pi_\lambda(f^\lambda) $ where for any $ g \in L^1\cap L^2(\R^{2n})$ the operator
$$ \pi_\lambda(g) =\int_{\R^{2n}} g(x,y)\, \pi_\lambda(x,y)\, dx\, dy$$
is known as the Weyl transform of $ g.$ It is known that  $ \pi_\lambda : g \rightarrow \pi_\lambda(g) $ is a constant multiple of a unitary operator from $ L^2(\R^{2n}) $ onto the space $ \mathcal{S}_2 $ consisting of all Hilbert-Schmidt operators on $ L^2(\R^n).$ The convolution between two functions $ f $ and $ g $ on $\He^n$ is defined by
$$ f \ast g(x,y,t) = \int_{\He^n} f((x,y,t)(u,v,s)^{-1})\, g(u,v,s) du\, dv\, ds.$$
This gives rise to a family of convolution structures on $ L^1(\R^{2n}) $ defined by the relation $ (f \ast g)^\lambda(x,y) = f^\lambda \ast_\lambda g^\lambda(x,y) .$ These are known as $\lambda$-twisted convolutions and  we make use of the relation $ \pi_\lambda(f^\lambda \ast_\lambda g^\lambda) = \pi_\lambda(f^\lambda) \pi_\lambda(g^\lambda) $ in  this article.

\subsection{Hermite and special Hermite semigroups} The spectral decomposition of the Hermite operator $ H(\lambda) = -\Delta+\lambda^2 |x|^2 $ is given by
$$ H(\lambda) = \sum_{k=0}^\infty (2k+n)|\lambda|\, P_k(\lambda) $$
where the spectral projection $ P_k(\lambda) $ onto the $k$-th eigenspacce associated to the eigenvalue $ (2k+n)|\lambda| $ is given by
$$ P_k(\lambda)f = \sum_{|\alpha|=k}\,(f, \Phi_\alpha^\lambda)\, \Phi_\alpha^\lambda.$$
Here $ \Phi_\lambda^\lambda, \, \alpha \in \mathbb N^n $ are the normraised, scaled Hermite functions which are eigenfunctions of $ H(\lambda) $ with eigenvalues $ (2|\alpha|+n)|\lambda|.$ It is well known that they form an orthonormal basis for $ L^2(\R^n).$ The Hermite semigroup is defined on $ L^p(\R^n) $ by
$$ e^{-tH(\lambda)}f = \sum_{k=0}^\infty \, e^{-t(2k+n)|\lambda|}\, P_k(\lambda)f.$$
 The semigroup $ e^{-tH(\lambda)}$ is an integral operator with an explicit kernel from which it follows that it is a contraction on $ L^p(\R^n) $ for  any $ 1 \leq p \leq \infty.$\\
 
 The special Hermite operator denoted by $ L_\lambda $ is the operator on $ \R^{2n} $ given by
 $$ L_\lambda = -\Delta_{\R^{2n}}+\frac{\lambda^2}{4}(|x|^2+|u|^2)+i\lambda \sum_{j=1}^n \left( x_j \frac{\partial}{\partial u_j}- u_j \frac{\partial}{\partial x_j} \right).$$
 There is an orthonormal basis for $ L^2(\R^{2n}) $ consisting of joint eigenfunctions of $ L_\lambda$  and $ H(\lambda)$ on $\R^{2n}.$ These are defined by 
 $$ \Phi_{\alpha \beta}^\lambda(x,u) = (2\pi)^{-n/2} |\lambda|^{n/2}\, (\pi_\lambda(x,u)\Phi_\alpha^\lambda, \Phi_\beta^\lambda).$$
\begin{equation}\label{spl-eigen}
 L_\lambda  \Phi_{\alpha \beta}^\lambda = (2|\beta|+n)|\lambda| \,\, \Phi_{\alpha \beta}^\lambda,\,\,\, H(\lambda)\Phi_{\alpha \beta}^\lambda = (2|\alpha|+2|\beta|+2n)|\lambda|\, \, \Phi_{\alpha \beta}^\lambda.
 \end{equation}
 The special Hermite functions satisfy another interesting orthogonality relation under twisted convolution (see \cite{T2}):
 \begin{equation}\label{spl-ortho}
 \Phi_{\alpha \beta}^\lambda \ast_\lambda \Phi_{\mu \nu}^\lambda = \delta_{\beta \mu} \, \Phi_{\alpha \nu}^\lambda.
 \end{equation}
  The spectrum of $ L_\lambda $ consists of $ (2k+n)|\lambda|, k \in \mathbb N $ but in view of the above
 the eigenspaces are infinite dimensional. However, the spectral decomposition of $ L_\lambda $ has a compact form in terms of Laguerre functions $ \varphi_k^\lambda(x,u). $
\\
 
 The $\lambda$-twisted convolution between two  functions $ f $ and $ g $ on $ \R^{2n}$ defined earlier is explicitly given by the integral
 \begin{equation}\label{twist-con} f \ast_\lambda g(x,u) = \int_{\R^{2n}} f(a,b)\, g(x-a,u-b)\, e^{-i\frac{\lambda}{2}( u\cdot a- x \cdot b)}\, da\, db.
 \end{equation}
 With this notation, the spectral decomposition of $ L_\lambda$ is given explicitly by
 $$ L_\lambda f(x,u) = (2\pi)^{-n} |\lambda|^n\,  \sum_{\alpha \in \mathbb N^n} (2|\alpha|+n)|\lambda|\,\,\, f \ast_\lambda \Phi_{\alpha \alpha}^\lambda(x,u)$$
 which can be further simplified. It is known that
 $$ \sum_{|\alpha|=k}  \Phi_{\alpha \alpha}^\lambda(x,u) = \varphi_{k,\lambda}^{n-1}(x,u)  = L_k^{n-1}(\frac{1}{2}|\lambda|(|x|^2+|u|^2))\, e^{-\frac{1}{4}|\lambda|(|x|^2|+|u|^2|)}$$
 where $ L_k^{n-1}(t), t \geq 0 $  stand for  the Laguerre polynomials of type $(n-1).$
 In terms of these functions the spectral decomposition of $ L_\lambda $ can be written as 
 $$ L_\lambda f(x,u) = (2\pi)^{-n} |\lambda|^n\,  \sum_{k=0}^\infty (2k+n)|\lambda|\, f \ast_\lambda \varphi_{k,\lambda}^{n-1}(x,u).$$
 Consequently, the special Hermite semigroup generated by $ L_\lambda $ is given by the expansion
 $$ e^{-tL_\lambda} f(x,u) = (2\pi)^{-n} |\lambda|^n\,  \sum_{k=0}^\infty e^{-t (2k+n)|\lambda| }\, f \ast_\lambda \varphi_{k,\lambda}^{n-1}(x,u).$$
 From the above definition it is clear that $e^{-tL_\lambda} f = f \ast_\lambda p_t^\lambda $ where the kernel $ p_t^\lambda $ is given by the expansion
 $$ (2\pi)^{-n} |\lambda|^n\,  \sum_{\alpha \in \mathbb N^n} e^{-t(2|\alpha|+n)|\lambda|} \, \Phi_{\alpha \alpha}^\lambda(x,u)  = (2\pi)^{-n} |\lambda|^n\, \sum_{k=0}^\infty e^{-t(2k+n)|\lambda|}\, \varphi_{k,\lambda}^{n-1}(x,u). $$
 The above series can be summed and $ p_t^\lambda $ is given explicitly by the formula
 $$ p_t^\lambda(x,u) = (4\pi)^{-n}  \lambda^n (\sinh t\lambda)^{-n} e^{-\frac{1}{4}\lambda (\coth t\lambda)(|x|^2+|u|^2)}.$$

  \subsection{Connection between Hermite and special Hermite semigroups}  The relation between the Hermite semigroup $ e^{-tH(\lambda)} $ on $ L^2(\R^n) $ and the special Hermite semigroup $ e^{-tL_\lambda} $ is given by  the Weyl transform: $ e^{-tH(\lambda)} = \pi_\lambda(p_t^\lambda).$ This is a well known relation but there is yet another  curious relation between Hermite and  special Hermite semigroups both acting  $ L^2(\R^{2n})$ which plays an important role in the proof of our main theorem.  In the next proposition $ H(\lambda) $ stands for the Hermite operator on $ \R^{2n}.$
 
 \begin{prop}\label{herm-spl-herm} For any $ t >0$ we have $ e^{-tH(\lambda)}f = p_t^\lambda \ast_\lambda f \ast_\lambda p_t^\lambda $ for all $ f \in L^2(\R^{2n}).$
 \end{prop}
 \begin{proof} As the special Hermite functions $ \Phi_{\alpha \beta}^\lambda $ form an orthonormal basis for $L^2(\R^{2n}) $ it enough to check that
 $$ e^{-tH(\lambda)} \Phi_{\alpha \beta}^\lambda = p_t^\lambda \ast_\lambda \Phi_{\alpha \beta}^\lambda \ast_\lambda p_t^\lambda $$ for all $ \alpha$ and $ \beta.$ Since $ p_t^\lambda $ is given by the exansion
 $$ p_t^\lambda(x,u) = \sum_{\mu \in \mathbb N^n} e^{-t(2|\mu|+n)|\lambda|} \, \Phi_{\mu \mu}^\lambda(x,u),$$ 
 using  the orthogonality relations \eqref{spl-ortho} we get 
 $$p_t^\lambda \ast_\lambda \Phi_{\alpha \beta}^\lambda  = e^{-t(2|\alpha|+n)|\lambda|} \, \Phi_{\alpha \beta}^\lambda,\,\,\,  \Phi_{\alpha \beta}^\lambda \ast_\lambda p_t^\lambda  = e^{-t(2|\beta|+n)|\lambda|} \, \Phi_{\alpha \beta}^\lambda.$$
 Consequently, 
 $$p_t^\lambda \ast_\lambda \Phi_{\alpha \beta}^\lambda \ast_\lambda p_t^\lambda  = e^{-t(2|\alpha|+2|\beta|+2n)|\lambda|}\,  \Phi_{\alpha \beta}^\lambda.$$
 As $ \Phi_{\alpha \beta}^\lambda $ are eigenfunctions of $ H(\lambda) $ with eigenvalues $ (2|\alpha|+2|\beta|+2n)|\lambda| $ we also have 
 $$ e^{-tH(\lambda)} \, \Phi_{\alpha \beta}^\lambda = e^{-t(2|\alpha|+2|\beta|+2n)|\lambda|}\,  \Phi_{\alpha \beta}^\lambda.$$
 This completes the proof of the proposition.
 \end{proof}

 The above result allows us to extend the domain of definition of  the Hermite semigroup  from $ L^2(\R^{2n}) $ to a bigger space. Let  $ L^2(\R^{2n},\, p_t^\lambda) $ stand for the space of all tempered distributions $ f $ for which $ p_t^\lambda  \ast_\lambda f \in L^2(\R^{2n}).$ For such functions 
 $$ e^{-tH(\lambda)}f = p_t^\lambda \ast_\lambda f \ast_\lambda p_t^\lambda $$ 
 is well defined and extends to $ \C^{2n} $ as an entire function. 

\subsection{On the image of $L^2(\R^{2n}) $ under  Hermite and special Hermite semigroups}  For any $ f \in L^2(\R^{2n})$ the functions $ e^{-tH(\lambda)}f $ and $ e^{-tL_\lambda} $ both extend to $ \C^{2n} $ as entire functions. Thus the image of $ L^2(\R^{2n}) $ under either of these semigroups consist of entire functions which are square integrable with respect to certain explicit positive weight functions. For the case of the Hermite semigroup the weight function is given by
\begin{equation}\label{herm-weight} U_t^\lambda(z,w) = c_n (\sinh 4t \lambda)^{-2n} \, e^{- \lambda(\coth 2t \lambda) (|y|^2+|v|^2)+ \lambda (\tanh 2t\lambda)(|x|^2+|u|^2)} .
\end{equation}
Thus for any $ f \in L^2(\R^{2n})$ the function $ F(z,w) = e^{-tH(\lambda)}f(z,w) $ satisfies the identity
\begin{equation}\label{herm-berg} \int_{\C^n }\int_{\C^n} |F(z,w)|^2\, U_t^\lambda(z,w)\, dz\, dw = = c_{n,\lambda} \int_{\R^{2n}} |f(x,u)|^2\, dx\, du.
\end{equation}
We have stated the above result for $ H(\lambda) $ on $ \R^{2n}$ but the result is true in any dimension. This was first proved by Byun \cite{B} but there are many other proofs available in the literature, see \cite{T4}.  The  proof given in \cite{T4} relies on the so called Gutzmer's formula for the Hermite expansions which also plays an important role in getting a necessary condition on $\varphi$ so that $ S_\varphi^\lambda $ is bounded on the twisted Fock space.\\

In order to state Gutzmer's formula we need to set up some notation.  We make use of several facts about the Schr\"odinger representation $ \pi_\lambda $ of the Heisenberg group $\He^n.$  The operator $ \pi_\lambda(x,u) $ can also be defined for complex values of the arguments. Thus for $ z, w \in \C^n,\, \pi_\lambda(z,w) $ can be defined on $ L^2(\R^n) $ as densely defined unbounded operators. If we let $ U(n) $ denote the group of $ n \times n $ unitary matrices, then for any $ \sigma \in U(n),$ the map $ (z,t) \rightarrow (\sigma z,t)$ is an automorphism of $ \He^n.$  The Laguerre function $ \varphi_{k,\lambda}^{n-1}(x,u) $ can also be extended to $ \C^n \times \C^n$ as an entire function. With these notations we have 

\begin{thm}[Gutzmer]\label{Gutz}
Let $ f \in L^2(\R^{n}) $ has a holomorphic extension $ F $ to $ \C^{n}.$ Then  we have
$$ \int_{U(n)}   \int_{R^{n}}  | \pi_\lambda(i\sigma(y,v))F(\xi)|^2   \, d\xi\, d\sigma 
 = c_n  \sum_{k=0}^\infty  \frac{k!(n-1)!}{(k+n-1)!} \varphi_{k,\lambda}^{n-1}(2iy,2iv)\,  \| P_k(\lambda)f \|_2^2\,.$$
\end{thm}

We refer to \cite{T4} for a proof of this result.  This formula along with the following lemma leads to a proof of the identity \eqref{herm-berg} in the case of $ e^{-tH(\lambda)} $ on $ L^2(\R^n).$ 

\begin{lem}\label{lemma} For any $ t > 0$ we have
$$  \frac{k!(n-1)!}{(k+n-1)!} \,\int_{\R^{2n}} p_{2t}^\lambda(2y,2v)\,\varphi_{k,\lambda}^{n-1}(2iy,2iv)\,\, dy\, dv = c_\lambda \, e^{ 2t(2k+n)|\lambda|}.$$
\end{lem}

\noindent It is interesting to compare the above  with the following identity which is easy to verify:
$$  \frac{k!(n-1)!}{(k+n-1)!} \,\int_{\R^{2n}} p_{2t}^\lambda(2y,2v)\,\varphi_{k,\lambda}^{n-1}(2y,2v)\,\, dy\, dv = c_\lambda \, e^{ -2t(2k+n)|\lambda|}$$
In fact, this is an immediate consequence  of the exapnsion of $ p_t^\lambda $ in terms of Laguerre functions $ \varphi_{k,\lambda}^{n-1}$ and their orthogonality in $ L^2(\R^{2n}).$\\

In the case of the special Hermite semigroup we have a similar result. For $ f \in L^2(\R^{2n})$ the function $ F(x,u) = e^{-tL_\lambda}f(x,u) $ has a holomorphic extension to $ \C^{2n}$ and we have the identity
\begin{equation}\label{spl-herm-berg}
\int_{\C^{n}} \int_{\C^{n}} |F(z,w)|^2 \, W_t^\lambda(z,w)  \, dz\,dw = c_{\lambda,n} \int_{\R^{2n}} |f(x,u)|^2 \, dx\, du
\end{equation}
where the weight function $ W_t^\lambda(z,w) $ is given by
$$ W_t^\lambda(z,w)= e^{\lambda \Im (z \cdot \bar{w})}\, p_{2t}^\lambda(2y,2v),\,\,\, z = x+iy, w=u+iv.$$
Thus the image of $ L^2(\R^{2n})$ under the special Hermite semigroup is precisely the weighted Bergman space consisting of entire functions on $ \C^{2n}$ for which \eqref{spl-herm-berg} is finite.
This was first proved in \cite{KTX} but there is   another proof using Gutzmer's formula for special Hermite expansion. We refer the reader  to the paper \cite{T3} for details. We conclude this subsection with the following remarks.

\begin{rem}  Note that in view of  Proposition \ref{herm-spl-herm} and \eqref{spl-herm-berg}, for any $ f \in L^2(\R^{2n}, p_t^\lambda)$ we have the identity
$$ \int_{\C^{2n}} | e^{-tH(\lambda)}f(z,w)|^2 \, W_t^\lambda(z,w)\, dz\, dw = \int_{\R^{2n}} | p_t^\lambda \ast_\lambda f(x,u)|^2 dx du.$$
As both $ H(\lambda) $ and $ p_t^\lambda$ are even functions of $ \lambda,$ recalling the definition of $ W_t^\lambda $ and using the fact that $  p_t^\lambda \ast_{-\lambda} f = f \ast_\lambda p_t^\lambda $ we have
$$ \int_{\C^{2n}} | e^{-tH(\lambda)}f(z,w)|^2\, e^{-\lambda \Im (z \cdot \bar{w})}\, p_{2t}^\lambda(2y,2v) dz\, dw = \int_{\R^{2n}} | f \ast_\lambda p_t^\lambda(x,u)|^2 dx du.$$
Using the identity \eqref{spl-herm-berg} once more we can rewrite the right hand side in terms of $ f \ast_\lambda p_{2t}^\lambda $ leading to the interesting relation
$$ \int_{\C^{2n}} | e^{-tH(\lambda)}f(z,w)|^2\,W_t^{-\lambda}(z,w) dz\, dw =  \int_{\C^{2n}} | e^{-2tL_\lambda}f(z,w)|^2\, W_{2t}^\lambda(z,w) dz\, dw .$$
\end{rem}

\begin{rem}
Let $ \mathcal{S}_2(\Gamma_\lambda) $ be the image of $ L^2(\R^{2n}, p_t^\lambda) $ under the Weyl transform. In other words, a densely defined linear operator $ T $ on $ L^2(\R^n) $ belongs to $ \mathcal{S}_2(\Gamma_\lambda) $  if $ e^{-\frac{1}{2}H(\lambda)} T $ is Hilbert-Schmidt. Note that $ B(L^2(\R^n)) \subset \mathcal{S}_2(\Gamma_\lambda)  $ and $ G_\lambda $ initially defined on $B(L^2(\R^n)) $ has a natural extension to $ \mathcal{S}_2(\Gamma_\lambda) .$ Thus if $T =\pi_\lambda(f)$ we have
$$ G_\lambda(T)(z,w) = p_1^\lambda(z,w)^{-1} e^{-\frac{1}{2}H(\lambda)} f(z,w).$$ 
We denote the image of  $ \mathcal{S}_2(\Gamma_\lambda)$ under $G_\lambda$ by $ \mathcal{F}_0^\lambda(\C^{2n}).$
\end{rem}

\begin{rem} When $ T =  \pi_\lambda(f) \in B(L^2(\R^n))$ we have the estimate
$$ \int_{\C^{2n}} | e^{-\frac{1}{2}H(\lambda)}f(z,w)|^2 \,W_t^\lambda(z,w)\, dz dw = \| T e^{-\frac{1}{2}H(\lambda)}\|_{HS}^2  \leq C \| T\|_{op}^2.$$
In particular, when $ f \in L^p(\R^{2n}), 1 \leq p \leq 2, \| \pi_\lambda(f) \|_{op} \leq \|f \|_p $  and hence we  have
$$ \int_{\C^{2n}} | e^{-tH(\lambda)}f(z,w)|^2 \,W_t^\lambda(z,w)\, dz dw  \leq C \| f\|_p^2 .$$
\end{rem}

\subsection{Necessary and sufficient conditions for $ S_\varphi^\lambda$ to be bounded}
 From the  above discussions, it is clear that a necessary and sufficient condition for  $ S_\varphi^\lambda $ to be  bounded on $ \mathcal{F}^\lambda(\C^{2n}) $ is that $ \varphi(z,w) = p_1^\lambda(z,w)^{-1} e^{-\frac{1}{2}H(\lambda)}f(z,w) $ with $ \pi_\lambda(f) \in B(L^2(\R^n)).$ This in turn is equivalent to the statement that the operator $ T_f: g \rightarrow g \ast_\lambda f $ initially defined on $ \mathcal{S}(\R^{2n}) $ extends to $ L^2(\R^{2n}) $ as a bounded operator if and only if $ \pi_\lambda(f) \in B(L^2(\R^n)).$ 
By using Gutzmer's formula for  Hermite expansions, we can  get a necessary condition on $ \varphi.$ In what follows we write  $ U_t^\lambda(\xi, \eta)$ in place of $ U_t^\lambda(z,w) $ where $ (z,w) = \zeta = \xi+i\eta .$

\begin{prop} For $ f \in L^2(\R^{2n}, p_{1/2}^\lambda)$ let $ \varphi = G_\lambda(\pi_\lambda(f))$ and define $ F = p_1^\lambda \, \varphi.$ Then for any $ 0 <t <1/2,$ we have the identity
$$ \int_{\R^{4n}} |F(\xi+i\eta)|^2 \, U_t^\lambda(\xi,\eta) \, d\xi\, d\eta = c_n  \sum_{k=0}^\infty \, \sum_{|\alpha|+|\beta|=k} | \left(\pi_\lambda(f)\Phi_\alpha^\lambda, \Phi_\beta^\lambda \right) |^2 \,  e^{-(1-2t)(2k+2n)|\lambda|}.$$
\end{prop}
\begin{proof}
Suppose $ \varphi(z,w) = B_\lambda (p_{1/2}^\lambda \ast_\lambda f)(z,w) $  where $ f \in L^2(\R^{2n}, p_{1/2}^\lambda).$  Then we have
$$ F(z,w) = \varphi(z,w)\, p_1^\lambda(z,w) = e^{-\frac{1}{2}H(\lambda)}f(z,w) .$$
Thus $ F $ is the holomorphic extension of  the  function $ e^{-\frac{1}{2}H(\lambda)}f \in L^2(\R^{2n}).$ By Gutzmer's formula  stated in Theorem \ref{Gutz} we obtain
$$ \int_{U(2n)}   \int_{R^{2n}}  | \pi_\lambda(i\sigma(\eta^\prime, \eta))F(\xi)|^2   \, d\xi\, d\sigma $$
$$  = c_n  \sum_{k=0}^\infty  e^{-(2k+2n)|\lambda|}\,  \frac{k!(2n-1)!}{(k+2n-1)!} \varphi_{k,\lambda}^{2n-1}(2i(\eta^\prime,\eta))\,  \| P_k(\lambda)f \|_2^2\,.$$
As $ H(\lambda)$ is a differential operator with real coefficients, we note that 
$$ H(\lambda) \overline{\Phi}_{\alpha \beta}^\lambda = \left(2|\alpha|+2|\beta|+2n\right)|\lambda|\,\overline{\Phi}_{\alpha \beta}^\lambda .$$
Hence $ P_k(\lambda)f$ which is the projection of $ f $ into the $k$-the eigenspace corresponding to the eigenvalue $(2k+2n)|\lambda|$ is given by two expressions, viz.
$$ P_k(\lambda)f = \sum_{|\alpha|+|\beta|=k} (f, \Phi_{\alpha \beta}^\lambda)\, \Phi_{\alpha \beta}^\lambda = \sum_{|\alpha|+|\beta|=k} (f, \overline{\Phi}_{\alpha \beta}^\lambda)\, \overline{\Phi}_{\alpha \beta}^\lambda. $$
We use the second expression to calculate $ \| P_k(\lambda)f\|_2^2.$ Recalling the definition of the special Hermite functions we see that
$$ \big(f, \overline{\Phi}_{\alpha \beta}^\lambda \big) = (2\pi)^{-n/2}\int_{\R^{2n}} f(x,u) \big(\pi_\lambda(x,u)\Phi_\alpha^\lambda, \Phi_\beta^\lambda \big) dx\,du = (2\pi)^{-n/2}\big(\pi_\lambda(f)\Phi_\alpha^\lambda, \Phi_\beta^\lambda \big) $$
and hence we have
$$ \| P_k(\lambda)f\|_2^2 =  (2\pi)^{-n} \sum_{|\alpha|+|\beta|=k} | \big(\pi_\lambda(f)\Phi_\alpha^\lambda, \Phi_\beta^\lambda \big) |^2.$$
Therefore, we have proved the identity
$$ \int_{U(2n)}   \int_{R^{2n}}  | \pi_\lambda(i\sigma(\eta^\prime,\eta))F(\xi)|^2   \, d\xi\, d\sigma $$
$$ = c_n  \sum_{k=0}^\infty  e^{-(2k+2n)|\lambda|}\,  \frac{k!(2n-1)!}{(k+2n-1)!} \varphi_{k,\lambda}^{2n-1}(2i(\eta^\prime,\eta))\,  \sum_{|\alpha|+|\beta|=k} | \big(\pi_\lambda(f)\Phi_\alpha^\lambda, \Phi_\beta^\lambda \big) |^2.$$
As $ p_{2t}^\lambda(2(\eta^\prime,\eta)) $ is invariant under the action of $ U(2n) $ the above leads to the identity
$$ \int_{R^{4n}}  \left( \int_{R^{2n}}  | \pi_\lambda(i(\eta^\prime,\eta))F(\xi)|^2   \, d\xi\, \right) p_{2t}^\lambda(2(\eta^\prime,\eta))d\eta d\eta^\prime $$
$$= \int_{\R^{2n}} \left( \int_{U(2n)}   \int_{R^{2n}}  | \pi_\lambda(i\sigma(\eta^\prime,\eta))F(\xi)|^2   \, d\xi\, d\sigma\right)  p_{2t}^\lambda(2(\eta^\prime,\eta))d\eta d\eta^\prime.$$
The integral on the left can be simplified and easily seen to be 
$$ \int_{\R^{4n}} |F(\xi+i\eta)|^2 \, U_t^\lambda(\xi,\eta) \, d\xi\, d\eta$$
whereas in view of the identity proved above the integral on the right is given by 
$$  \sum_{k=0}^\infty  \, \sum_{|\alpha|+|\beta|=k} | \big(\pi_\lambda(f)\Phi_\alpha^\lambda, \Phi_\beta^\lambda \big) |^2 \, c_{k,\lambda}\, e^{-(2k+2n)|\lambda|} $$
where we have written $ c_{k,\lambda} $  for the constant defined by the following integral:
$$  c_{k,\lambda} = \frac{k!(2n-1)!}{(k+2n-1)!}  \left(\int_{\R^{4n}}\varphi_{k,\lambda}^{2n-1}(2i(\eta^\prime,\eta))\, p_{2t}^\lambda(2(\eta^\prime,\eta))d\eta d\eta^\prime\right).$$
We now make use of the identity stated in Lemma \ref{lemma} to conclude that
$$ \int_{\R^{4n}} |F(\xi+i\eta)|^2 \, U_t^\lambda(\xi,\eta) \, d\xi\, d\eta = c_n  \sum_{k=0}^\infty \,\sum_{|\alpha|+|\beta|=k} | \big(\pi_\lambda(f)\Phi_\alpha^\lambda, \Phi_\beta^\lambda \big)  |^2 \,e^{-(1-2t)(2k+2n)|\lambda|}. $$
This completes the proof of the proposition.
\end{proof}

Once we have the above proposition at our disposal,  it is  easy to get a necessary condition and a slightly stronger sufficient condition on $ \varphi$ so that $ S_\varphi^\lambda $ is bounded. With $ F(z,w) = p_1^\lambda(z,w)\varphi(z,w) $ we can easily check that $ F $ is square integrable with respect to $ U_t^\lambda(\xi,\eta)$ if and only if $ \varphi(\xi+i\eta)=\varphi(z,w) $ is square integrable with respect to the weight function
$$ U_t^\lambda(\xi,\eta) \, |p_1^\lambda(\xi+i\eta)|^2 = U_t^\lambda(\xi,\eta)\, e^{-\frac{\lambda}{2}(\coth \lambda)|\xi|^2}\, e^{\frac{\lambda}{2}(\coth \lambda)|\eta|^2}.$$
After simplification the above reduces to the weight function
 $ w_t^\lambda(\xi,\eta) $ which appears in Theorem 1.3 and defined in \eqref{weight}.  Therefore, Theorem 1.3  will follow once we prove the following result. Let us recall that  $ \mathcal{F}_0^\lambda(\C^{2n})$ stands for  the image of  $ \mathcal{S}_2(\Gamma_\lambda)$ under $G_\lambda.$ \\

\begin{thm} Given $ \varphi \in \mathcal{F}_0^\lambda(\C^{2n}),$ let $ F = p_1^\lambda\, \varphi.$ A necessary condition for $ S_\varphi^\lambda $ to be bounded on $\mathcal{F}^\lambda(\C^{2n})$ is that
$$ \int_{\R^{4n}} |F(\xi+i\eta)|^2 \, U_t^\lambda(\xi,\eta) \, d\xi\, d\eta < \infty $$
for all $ 0 < t < 1/2.$ Conversely, if the above condition holds for $ t = 1/2,$ then $ S_\varphi^\lambda $ is bounded.
\end{thm} 
\begin{proof}

Since $ \varphi  = G_\lambda(\pi_\lambda(f))$  for some $ f \in L^2(\R^{2n}, p_{1/2}^\lambda) $ it follows that $ S_\varphi^\lambda $ is bounded on $ \mathcal{F}^\lambda(\C^{2n}) $ if and only if $ \pi_\lambda(f) $ is bounded on $ L^2(\R^n).$ When this happens
$$ \sum_{|\alpha|+|\beta|=k} | \big(\pi_\lambda(f)\Phi_\alpha^\lambda, \Phi_\beta^\lambda \big) |^2 \leq C \,  \frac{(k+2n-1)!}{k!(2n-1)!} $$
and consequently from the proposition proved above we get the estimate
$$ \int_{\R^{4n}} |F(\xi+i\eta)|^2 \, U_t^\lambda(\xi,\eta) \, d\xi\, d\eta = C_n  \sum_{k=0}^\infty \,    \frac{(k+2n-1)!}{k!(2n-1)!} \,  e^{-(1-2t)(2k+2n)|\lambda|} $$
which is finite as long as $ 0 < t < 1/2.$ Conversely, if the condition holds for $  t = 1/2$ we have
$$ \sum_{k=0}^\infty \, \sum_{|\alpha|+|\beta|=k} | \left(\pi_\lambda(f)\Phi_\alpha^\lambda, \Phi_\beta^\lambda \right) |^2  < \infty.$$
But  this simply means that $ \pi_\lambda(f) $ is a Hillbert-Schmidt operator since
$$ \| \pi_\lambda(f) \|_{HS}^2 = \sum_{\alpha \in \mathbb N^n} \| \pi_\lambda(f)\Phi_\alpha^\lambda \|_2^2=  \sum_{\alpha, \beta \in \mathbb N^n} | \left(\pi_\lambda(f)\Phi_\alpha^\lambda, \Phi_\beta^\lambda \right) |^2.$$
This completes the proof of the theorem and hence Theorem 1.3 is also proved.
\end{proof}

\medskip


\section*{Acknowledgements}
This work was carried out when the author was visiting Harish-Chandra Research Institute as Infosys Visiting Professor. He thanks HRI for the facilities and warm hospitality. He  was also partially supported by INSA.
\\

\providecommand{\bysame}{\leavevmode\hbox to3em{\hrulefill}\thinspace}
\providecommand{\MR}{\relax\ifhmode\unskip\space\fi MR }
\providecommand{\MRhref}[2]{%
  \href{http://www.ams.org/mathscinet-getitem?mr=#1}{#2}
}
\providecommand{\href}[2]{#2}


\begin{thebibliography}{10}

\bibitem{B} D.-W. Byun, \emph{Inversions of Hermite semigroup,} Proc. Amer. Math. Soc. \textbf{118} (1993), 437–445..

\bibitem{BD} S. Bais and V. N. Dogga, \emph{$\mathcal{L}$-invariant radial singular integral operators on Fock spaces}, {J. Pseudo-Differ. Oper. Appl.} \textbf{14} (2023), no. 1, Paper No. 11, 35 pp.

\bibitem{CHLS}
G. Cao, L.~He, J.~Li, and M. Shen, \emph{Boundedness criterion for integral operators on the fractional {F}ock-{S}obolev spaces}, Math. Z.
  \textbf{301} (2022), no.~4, 3671--3693. \MR{4449725}

\bibitem{CLSWY}
G. Cao, J.~Li, M. Shen, B. D. Wick, and L. Yan, \emph{A boundedness criterion for singular integral operators of convolution type on
  the {F}ock space}, Adv. Math. \textbf{363} (2020), 107001, 33. \MR{4053517}

\bibitem{GBF}
G.~B. Folland, \emph{Harmonic analysis in phase space}, Annals of Mathematics Studies, vol. 122, Princeton University Press, Princeton, NJ, 1989. \MR{983366}
  
\bibitem{GT} R. Garg and S. Thangavelu, \emph{Boundedness of certain linear operators on twisted Fock spaces}, Math. Z. (2024)  arXiv:2306.14188 (2023)
  

\bibitem{KTX} B. Kr\"{o}tz, S. Thangavelu, and Y. Xu, \emph{The heat kernel transform for the {H}eisenberg group}, J. Funct. Anal. \textbf{225} (2005),  no.~2, 301--336. \MR{2152501}


\bibitem{T1} S. Thangavelu, \emph{Lectures on {H}ermite and {L}aguerre expansions}, 
Mathematical Notes, vol.~42, Princeton University Press, Princeton, NJ, 1993, With a preface by Robert S. Strichartz. \MR{1215939}

\bibitem{T2}  S. Thangavelu,  \emph{An introduction to the uncertainty principle}, Progress in
  Mathematics, vol. 217, Birkh\"{a}user Boston, Inc., Boston, MA, 2004, Hardy's  theorem on Lie groups, With a foreword by Gerald B. Folland. \MR{2008480}
  
  \bibitem{T3} S. Thangavelu, \emph{Gutzmer's formula and Poisson integrals on the Heisenberg
group,} Pacific J. Math. \textbf{231} (2007),  no.1, 217-237.
  
  \bibitem{T4} S. Thangavelu, \emph{An analogue of Gutzmer's formula for  Hermite expansions}, Studia Math. \textbf{185} (2008), 279-290.

\bibitem{T5} S. Thangavelu,  \emph{Fourier multipliers and pseudo-differential operators on
  {F}ock-{S}obolev spaces}, arXiv: 2304.01087 (2023).
  
 \bibitem{T6}  S. Thangavelu, \emph{Algebras of entire functions and  representations of the twisted Heisenberg group,}  Indian J. Pure and Appl. Math. (2024).
  
 \bibitem{WW} B. D. Wick and  S. K. Wu, \emph{Integral operators on Fock-Sobolev spaces via multipliers on Gauss-Sobolev spaces}, {Integr. Equ. Oper. Theory} \textbf{94} (2022), no. 2, Paper No. 22, 24 pp.
  

\end{thebibliography}
\end{document}